\documentclass[11pt]{article} 

\usepackage[utf8]{inputenc} 
\usepackage{ulem}

\usepackage{geometry} 
\usepackage{graphicx} 

\usepackage{booktabs} 
\usepackage{array} 
\usepackage{paralist} 
\usepackage{verbatim} 
\usepackage{subfig} 
\usepackage{amssymb}
\usepackage{amsmath}
\usepackage{amsthm}

\usepackage{fancyhdr} 
\usepackage{sectsty}
\allsectionsfont{\sffamily\mdseries\upshape} 
\theoremstyle{plain}
\newtheorem{thm}{Theorem}[section]

\newtheorem{lem}[thm]{Lemma}

\usepackage[nottoc,notlof,notlot]{tocbibind} 
\usepackage[titles,subfigure]{tocloft} 

\newcommand{\NN}{\mathbb{N}}

\title{Berlekamp's Switching Game on Finite Projective and Affine Planes}
\author{James Martin}

\begin{document}
\maketitle

\begin{abstract}
We adapt Berlekamp's light bulb switching game to finite projective plans and finite affine planes, then find the worst arrangement of lit bulbs for planes of even and odd orders. The results are then extended from the planes to spaces of higher dimension.
\end{abstract}

\section{Berlekamp's Switching Game}

The original game is simple enough. 100 light bulbs arranged in 10 rows and 10 columns. Each column and row has a switch that toggles the states of their 10 bulbs. 

Before the game begins, each bulb is individually switched on or off, giving the {\sl initial configuration}. The game is played by flipping the row and column switches with the goal of ending up with fewer lit bulbs than the game started with. If this can be achieved, the initial configuration is said to be {\sl reducible}. A board is said to be {\sl reduced} if it has fewer bulbs lit than its initial configuration.

Naturally this raises the question, what is the worst case? In other words, what is the maximum number of lit bulbs which cannot be reduced? For a 10x10 game, the answer is 35. Generalizing this to $n \times n$ matrices, exact answers are known for $n=2,...,12$. Lower bounds are known for $13,...,20$. See \cite{MR2094708} and \cite{MR992740}. Theoretical bounds exist for $n>20$. Finding a new exact answer using known methods is computationally intensive given that there are $2^{n^2}$ possible starting configurations.

\section{The Game on Projective Planes}

What about other arrangements of bulbs and switches? First consider a finite projective plane $\mathbb{P}$ of order $p$.  $\mathbb{P}$ contains $p^2 + p + 1$ points and $p^2 + p + 1$ lines. Each point lies on $p+1$ lines, and each line contains $p+1$ points. Any two points lie on exactly one line. The simplest plane is when $p=2$, and is called the Fano plane. An account of projective planes can be found in reference \cite{MR0333959}.

Consider a bulb at each point and a switch for each line. A switch toggles the state of the $n+1$ bulbs on its associated line. The question of worst case is far easier to answer now.

\begin{lem}
\label{thm:projectivereducible}
For $p$ prime, the board $\mathbb{P}=PG(2,p)$ can be reduced if it has at least 2 bulbs lit.

\begin{proof}
If $p=2$ then the reduction is straightforward. To see this, notice that each line contains 3 bulbs, and any 2 bulbs lie on exactly one line. So choose a line with at least 2 lit bulbs and flip the switch. The line will have at most 1 lit bulb, and the board is reduced.

if $p>2$ then $p$ is odd. Note that $(1)$ every point lies on $p+1$ lines and $(2)$ any 2 points lie on exactly one line.
So given any 2 light bulbs, labeled $A$ and $B$, there is one line $l$ containing both. $A$ lies on $p$ other lines, call them $A_1,...,A_p$. Likewise, $B$ lies on $p$ other lines, call them $B_1,...,B_p$.

Flip the switches for lines $A_1,...,A_p$ and $B_1,...,B_p$. Both $A$ and $B$ were initially on, each was toggled $p$ times, turning them off.
Any point on line $l$ other than $A$ and $B$ was not toggled at all. Finally, any point not on line $l$ must lie on exactly one of $A_1,...,A_p$ and on exactly one of $B_1,...,B_p$, so it would be toggled twice and returned to its original state.

$A$ and $B$ were switched from on to off, no other bulbs changed state, thus the board is reduced.
\end{proof}
\end{lem}

\begin{lem}
\label{thm:projectiveparity}
For $p$ odd, the parity of the number of lit bulbs does not change by flipping a switch.
\begin{proof}
Flipping a switch toggles the state of $p+1$ bulbs, an even number. Thus the change in number of lit bulbs is a multiple of two. It follows that a board with odd parity will remain odd, and a board with even parity will remain even.
\end{proof}
\end{lem}

\begin{thm}
\label{cor:oddparity}
For $p$ odd, any board with an even number of lit bulbs can be reduced to 0 lit bulbs, and any board with an odd number of lit bulbs can be reduced to 1 lit bulb. In particular, the worst case scenario for a plane with odd order is one lit bulb.
\begin{proof}
This is immediate from \ref{thm:projectivereducible} and \ref{thm:projectiveparity}.
\end{proof}
\end{thm}

\section{Projective Planes of Even Order}

By Theorem \ref{cor:oddparity} the only planes which could have a worst case scenario of more than one bulb are planes of order $2^{n}$ for some integer $n>1$. Naive brute force computation has shown that the unique (up to isomorphism) plane of order 4 has a worst case scenario of 6 lit bulbs. In this case, the bulb configuration has the maximum possible number of bulbs lit such that no single line is trivially reducible. Turning on any unlit bulb results in a configuration that can be reduced with one switch flip.

We conjecture that such a maximal configuration will be irreducible for every projective plane of order $2^{n}$ for $n>1$. In which case, these configurations will always be the worst case scenario since lighting just one more bulb results in a trivially reducible configuration.

\section{The Game on Affine Planes}
The {\sl affine plane} $\mathbb{A}$ can be thought of as $\mathbb{P}$ with a single line and it's incident points removed. $\mathbb{A}$ contains $p^2$ points and $p^2 + p$ lines. Each point lies on $p+1$ lines, and each line contains $p$ points. Any two points lie on exactly one line. However, unlike in a projective plane, an affine plane contains parallel lines. That is, given a line $L$ and a point $P$ not on $L$, there is a line through $P$ which does not intersect $L$.

\begin{thm}
\label{thm:affineodds}
For affine plane $\mathbb{A}$ with odd order $p$, the board $\mathbb{A}$ can be reduced if it has at least one bulb lit.
\begin{proof}
Suppose $\mathbb{A}$ has at least one bulb illuminated.

Label one lit bulb $A$. Every point lies on $p^{n}+1$ lines, so enumerate the lines incident to $A$, calling them $A_1,...,A_{p^{n}+1}$.

Toggle lines $A_1$ through $A_{p^{n}}$, leaving line $A_{p^{n}+1}$ unswitched. $p^{n}$ is odd so $A$ was turned off, and the remaining bulbs of $A_{p^{n}+1}$ are unchanged, but every other bulb of $\mathbb{A}$ was toggled exactly once.

Line $A_1 \setminus \{A\}$ has exactly $p^{n}$ points, and for each of these points there is a line through that point which is parallel to $A_{p^{n}+1}$. Label these lines $B_1,...,B_{p^{n}}$. 

The axioms for affine planes imply that the parallel relation between lines is transitive. Since $B_1$ through $B_{p^{n}}$ are all parallel to $A_{p^{n}+1}$, they are all pairwise parallel. Hence switching lines $B_1,...,B_{p^{n}}$ will toggle every bulb once, except those on $A_{p^{n}+1}$.

So every bulb in $\mathbb{A}$ was toggled either 0 times or 2 times, except for bulb $A$, which was toggled on odd number of times. $A$ was switched off, no other bulbs changed state, therefore, the board is reduced.
\end{proof}
\end{thm}

This again leaves the affine planes of order $2^{n}$ for $n \in \NN$ as the only chance to find a nontrivial worst case. And we are not disappointed. It's easy to verify the worst case for $p=2$ is one lit bulb. The worst case for $p=4$ has been computed as 4 lit bulbs.

Just as with projective planes, the bulb configuration has the maximum possible number of bulbs lit such that no single line is trivially reducible. Turning on any unlit bulb results in a configuration that can be reduced with one switch flip.

We conjecture such a maximal configuration will be irreducible for every affine plane of order $2^{n}$ for $n>1$. In which case, these configurations will always be the worst case scenario since lighting just one more bulb results in a trivially reducible configuration.

\section{Higher Dimensional Spaces}
Theorems \ref{cor:oddparity} and \ref{thm:affineodds} were written with projective and affine planes in mind. This makes sense given that the original board game was essentially two-dimensional. However the properties used in the proofs are common to even dimensional projective and affine spaces, thus they also hold for projective and affine spaces of even dimension. 

\bibliographystyle{amsplain}	
\bibliography{Researchbib}		

\end{document}